\documentclass[a4paper,11pt,oneside,reqno]{amsart}
\usepackage{amsmath}
\usepackage{amsthm}
\usepackage{amssymb}
\usepackage{graphicx}
\usepackage{mathrsfs}
\usepackage{tikz}
\usetikzlibrary{matrix,arrows}

\theoremstyle{plain}

\newtheorem{teore}{Theorem}[section]
\newtheorem{defin}[teore]{Definition}
\newtheorem{lem}[teore]{Lemma}
\newtheorem{coro}[teore]{Corollary}
\newtheorem{propo}[teore]{Proposition}
\newtheorem{claim}{Claim}

\newtheorem{lemin}{Lemma}
\newtheorem*{claim*}{Claim}
\newtheorem*{prop*}{Proposition}
\newtheorem*{lem*}{Lemma}
\theoremstyle{remark}

\newtheorem{ejemplo}[teore]{{\sc Example}}
\newtheorem{notas}[teore]{{\sc Remark}}
\newcommand{\nrm}[1]{\|#1\|}

\newcommand{\prop}{\begin{propo}}
\newcommand{\fprop}{\end{propo}}
\newcommand{\cor}{\begin{coro}}
\newcommand{\fcor}{\end{coro}}

\newcommand{\defi}{\begin{defin}\rm}
\newcommand{\fdefi}{\end{defin}}
\newcommand{\eje}{\begin{ejemplo}}
\newcommand{\feje}{\end{ejemplo}}
\newcommand{\lema}{\begin{lem}}
\newcommand{\flema}{\end{lem}}
\newcommand{\teor}{\begin{teore}}
\newcommand{\fteor}{\end{teore}}
\newcommand{\nota}{\begin{notas}\rm}
\newcommand{\fnota}{ \end{notas}}
\newcommand{\clam}{\begin{claim}}
\newcommand{\fclam}{\end{claim}}
\newcommand{\clams}{\begin{claim*}}
\newcommand{\fclams}{\end{claim*}}

\newcommand{\lclam}{\begin{lclaim}}
\newcommand{\flclam}{\end{lclaim}}
\newcommand{\prucl}{\prue[Proof of Claim:]}
\newcommand{\fprucl}{\fprue}
\newcommand{\ben}{\begin{enumerate}}
\newcommand{\een}{\end{enumerate}}
\newcommand{\bit}{\begin{itemize}}
\newcommand{\eit}{\end{itemize}}

\newcommand{\mc}[1]{\mathcal{#1}}

\newcommand{\mr}[1]{\mathrm{#1}}

\newcommand{\casos}{\begin{itemize}}
\newcommand{\fcasos}{\end{itemize}\setcounter{cs}{1}}

\newcommand{\pe}{\preceq}

\newcommand{\conj}[2]{ \{ {#1}\,:\,{#2} \} }

\newcommand{\om}{\omega}

\newcommand{\ka}{\kappa}

\newcommand{\buit}{\emptyset}

\newcommand{\id}{\text{Id }}

\newcommand{\al}{\alpha}

\newcommand{\de}{\delta}
\newcommand{\De}{\Delta}
\newcommand{\la}{\lambda}

\newcommand{\N}{{\mathbb N}}

\renewcommand{\ker}{\mathrm{Ker}}
\newcommand{\rest}{\upharpoonright}

\newcommand{\con}{\subseteq}

\newcommand{\cones}{\varsubsetneq }

\newcommand{\prue}{\begin{proof}}
\newcommand{\fprue}{\end{proof}}
\makeindex
\begin{document}
\title{A Bourgain-Pisier construction for general Banach spaces}
\author{J. Lopez-Abad}

\address{Instituto de Ciencias Matematicas (ICMAT), CSIC-UAM-UC3M-UCM, Madrid, Spain.}
\email{abad@icmat.es}



\begin{abstract}
We prove that every Banach space, not necessarily separable, can be isometrically embedded into a $\mc
L_{\infty}$-space in a way that the corresponding quotient has  the Radon-Nikodym and the Schur properties.
As a consequence, we obtain  $\mc L_\infty$ spaces of arbitrary large densities with the Schur and the
Radon-Nikodym properties. This extents the  result by J. Bourgain and G. Pisier in \cite{BP} for separable
spaces.
\end{abstract}

\maketitle

\section{Introduction}
 The main question considered in this paper is the largeness of the class of $\mc L_\infty$ spaces
  in terms of embeddability.  Recall that a Banach space $X$ is called $\mc L_{\infty,\la}$ when for every
  finite dimensional subspace $F$ of $X$ there is a subspace $G$ of $X$ $\la$-isomorphic to $\ell_\infty^{\dim
  G}$ containing $F$. $\mc L_\infty$ just means $\mc L_{\infty,\la}$ for some $\la$.  There are two remarkable results
for the class of separable $\mc L_\infty$ spaces. The first  by J. Bourgain and G. Pisier in \cite{BP} states
that every separable Banach space $X$ can be isometrically embedded into a $\mc L_\infty$-space  $Y_X$ in
such a way that the corresponding quotient space $Y_X/X$ has the  Radon-Nikodym property (RNP) and the Schur
property. The second, more recent one, by D. Freeman, E. Odell and Th. Schlumprecht \cite{FOS}  tells that
every space with separable dual can be isomorphically embedded into a $\mc L_\infty$-space with separable
dual (therefore an $\ell_1$-predual). Both constructions are the natural extensions of  the work of J.
Bourgain and F. Delbaen \cite{BD} and Bourgain \cite{B}. There several other recent examples. Perhaps the
most impressive one is the $\mc L_\infty$-space by S. A. Argyros and R. G. Haydon \cite{AH} where every
operator is the sum of a multiple of the identity and a compact one.

In the non-separable context much less is known. Spaces of functions on a non-metrizable compactum, or
non-separable Gurarij spaces are non-separable $\mc L_\infty$-spaces.   There is a wide variety of structures
in the non-separable level for spaces in these two classes. Based on combinatorial axioms outside ZFC, there
are non-separable spaces in these two classes without uncountable biorthogonal systems or where every
operator is the sum of a multiple of the identity and an operator with separable range (see  \cite{LT} for
more information). On the other hand,   the separable structure of the known examples is too simple:  either
they are $c_0$-saturated, that is, every infinite dimensional subspace of it contains an isomorphic copy of
$c_0$, or universal for the separable spaces. So it is natural to ask if  there are examples of non-separable
$\mc L_\infty$-space without isomorphic copies of $c_0$, with the (RNP) or with the Schur property. Our main
result in  Theorem  \ref{maintheorem} is that the embedding  Theorem by Bourgain and Pisier reminds valid for
any density. In particular, by embedding $\ell_1(\ka)$ for an infinite cardinal number $\ka$, we obtain
examples of $\mc L_\infty$ spaces of arbitrary density with the Radon-Nikodym and the Schur properties.

For a given separable space $X$, the corresponding Bourgain-Pisier superspace $Y_X$ of it is  built  in such
a way that $Y_X$  and the quotient $Y_X/X$ are both the inductive limit of \emph{linear} systems
$(Z_n,j_n)_{n\in \N}$ of   a special type of isometrical embedding $j_n:Z_n\to Z_{n+1}$ ($\eta$-admissible
embeddings, see Definition \ref{ioioijo4trt}), and such that, in addition, the corresponding $Z_n$'s are
finite dimensional for the quotient space $Y_X/X$. The key fact to get  the Radon-Nikodym   and the Schur
properties of the quotient space $Y_X/X$ is the metric property of $\eta$-admissible embeddings exposed here
in Lemma \ref{kjdfklsdljfdw} and its consequence to inductive limits as above for finite dimensional spaces
(see \cite[Theorem 1.6]{BP}).

In contrast to the separable case, the main difficulty in the non-separable case is the construction the
appropriate inductive limit. Indeed, if $X$ is non-separable, then it is unlikely to find a nice linear
system having the space $Y_X$ as the corresponding limit. In general, every Banach space $X$ is naturally
represented as the inductive limit of its finite dimensional subspaces together with the corresponding
inclusions between them.  Our inductive system $((E_s)_{s\in I}, (j_{t,s})_{t \con s})$ to represent $Y_X$ is
also based on the inclusion relation over the index set $I$ consisting on all finite subsets of the density
of the  space $X$. This  provides a natural way to isometrically embed $X$ into $Y_X$.    In addition, our
inductive system is constructed in a way that its linear subsystems $((E_{s_n})_{n\in
\N},(j_{s_n,s_{n+1}})_{n\in \N})$ are Bourgain-Pisier linear systems as above. In other words, every
separable subspace $Z$ of $Y_X$ can be isometrically embedded into the separable Bourgain-Pisier extension
$Y_Z$.   So, having into account that  the Radon-Nikodym and Schur properties are separably determined, we
readily have that the quotient $Y_X/X$ has the desired properties.

To construct  the  spaces $E_s$ and the corresponding embeddings $j_{s,t}:E_s\to E_t$  we define first finite
linear systems  $((E_{t}^{(s)})_{t\con s}, (j_{u,t}^{(s)})_{u\prec t})$ of $\eta$-admissible embeddings
$j_{u,t}^{(s)}:E_{u}^{(s)}\to E_{t}^{(s)}$, where $\prec$ is a natural well ordering extending the inclusion
relation. Obviously, this raises   a  problem of coherence, since given  $s\con p\con q$ we will have defined
two ``$s$-extensions''  $E_{s}^{(p)}$ and $E_{s}^{(q)}$ of $E_\buit=X$  and therefore two isometric
embeddings $X\to E_s^{(p)}$ and $X\to E_s^{(q)}$. This is corrected by defining simultaneously an infinite
directed system $((E_{s}^{(p)})_{s\con p}, (k_{s}^{(p,q)})_{s\con p\con q})$ of $\eta$-admissible embeddings
making the appropriate diagrams commutative.

Finally, let us point out that nothing is known in how to skip the separability assumption in the
Freeman-Odell-Schlumprecht embedding Theorem, or, even more basic,  if a non-separable Bourgain-Delbaen
exists, i.e. a non-separable $\mc L_\infty$-space not containing isomorphic copies of $c_0$ or $\ell_1$.

The paper is organized as follows. The Section 2  is a survey of basic facts concerning the Kisliakov's
extension method and $\eta$-admissible, in particular we present a new extension fact concerning these
embeddings in Lemma \ref{i4jrijeerdthtyolo}. The last section is devoted to the proof of the Theorem
\ref{maintheorem}.

\section{Background and basic facts}
We use standard terminology in Banach space theory from the monographs \cite{AK} and \cite{LiTza}. The goal
of this section is to present the basic notions of $\eta$-admissible diagrams and $\eta$-admissible
embeddings introduced by Bourgain and Pisier.  To complete the information we give here, specially for some
proofs, we refer the reader to the original paper \cite{BP} or to the recent book by P. Dodos \cite{Do}.

Recall the Kisliakov's extension method \cite{K}: Given Banach spaces $S\con B$ and $E$ and an operator
$u:S\to E$ such that $\nrm{u}\le \eta\le 1$, let
\begin{align*}
N_u:= & \conj{(s,-u(s))\in B\times
E}{s\in S},\\
i_B:B& \to (B\oplus_1 E)/N_u \\
b &\mapsto i_B(b)=(b,0)+N_u \\
i_E:E &\to (B\oplus_1 E )/N_u \\
e&\mapsto i_E(e)=(0,e)+N_u
\end{align*}
Then the diagram $(K)$
\begin{equation*}
\begin{tikzpicture}[descr/.style={fill=white,inner sep=2pt}]
\matrix (m) [matrix of math nodes, row sep=3em, column sep=3em, text height=1.5ex, text depth=0.25ex]
{ B & (B\oplus_1 E)/N_u \\
S & E
\\
};

\path[->,font=\normalsize]

(m-1-1) edge node[above] {$i_B$} (m-1-2)

(m-2-1) edge node[below]  {$u$}  (m-2-2)

(m-2-2) edge node[right] {$i_E$} (m-1-2)

(m-2-1) edge node[right=26pt] {$(K)$}   (m-1-1)

 ;

\path[right hook->]

(m-2-1) edge     (m-1-1)

;

\end{tikzpicture}
\end{equation*}
is commutative, $i_E$ is an isometrical embedding, and $\nrm{i_B}\le 1$. This diagram has several categorical
properties such as minimality and uniqueness.
\defi
We say that a   diagram
\begin{equation*}
\begin{tikzpicture}[descr/.style={fill=white,inner sep=2pt}]
\matrix (m) [matrix of math nodes, row sep=3em, column sep=3em, text height=1.5ex, text depth=0.25ex]
{ B & E_1 \\
S & E
\\
};

\path[->,font=\normalsize]

(m-1-1) edge node[above] {$\bar{u}$} (m-1-2)

(m-2-1) edge node[below]  {$u$}  (m-2-2)

(m-2-2) edge node[right] {$j$} (m-1-2)

 ;

\path[right hook->]

(m-2-1) edge     (m-1-1)

;

\end{tikzpicture}
\end{equation*}
is a $\eta$-admissible diagram when there is an isometry  $T:(B\oplus_1 E)/N_{u}\to E_1$ such that $j=T\circ
i_E$ and $\bar{u}=T\circ i_B $.  The \emph{canonical} $\eta$-admissible diagram associated to the triple
$(S,B,u)$ is the Kisliakov's diagram $(K)$ above.

An isometrical embedding $j:E\to E_1$ is called \emph{$\eta$-admissible embedding} when there are $S\con B$,
$E_1$, $u:S\to E$, $\bar{u}: B\to E_1$ forming together with $j:E\to E_1$ an $\eta$-admissible diagram.
\fdefi
Observe that $\eta$-admissible diagrams are always commutative.

\defi\label{ioioijo4trt}\cite{BP}
A surjective operator  $\pi:E\to F$ is called a \emph{metric surjection} when the associated isomorphism $\bar{\pi}:E/\ker(\pi)\to F$ is an isometry.

\fdefi
The following are useful known characterizations, not difficult to prove.
\prop\label{jiejrijeijifjdd}
\begin{enumerate}
\item[(a)] Let $S\con B$, $E$ and $E_1$ be normed spaces, and $\eta\le 1$. A diagram $(\De)$
\begin{equation}
\begin{tikzpicture}[descr/.style={fill=white,inner sep=2pt}]
\matrix (m) [matrix of math nodes, row sep=3em, column sep=3em, text height=1.5ex, text depth=0.25ex]
{ B & E_1 \\
S & E
\\
};

\path[->,font=\normalsize]

(m-1-1) edge node[above] {$\bar{u} $} (m-1-2)

(m-2-1) edge node[below]  {$u$}  (m-2-2)

(m-2-2) edge node[right] {$j$} (m-1-2)

 ;

\path[right hook->]

(m-2-1) edge node[right=13pt] {$(\De)$}   (m-1-1)

;

\end{tikzpicture}
\end{equation}
is an \emph{$\eta$-admissible diagram} if and only if
\begin{enumerate}
\item[$(\al.1)$] $j$ is an isometry and $\nrm{u} \le \eta$,
\item[$(\al.2)$] $\pi:B\oplus_1 E \to E_1$ defined for $(b,e)\in B\times E$ by $\pi(b,e):=\bar{u}(b)+j(e)$ is
a metric surjection.
\item[$(\al.3)$] $\ker (\pi)=N_u=\conj{(s,-u(s))}{s\in S}$.
\end{enumerate}
\item[(b)] An isometrical embedding $j:E\to E_1$ is $\eta$-admissible iff there is some Banach space $B$ and a metric surjection
 $\pi:B\oplus_1 E \to E_1$ such that $\pi(0,e)=j(e)$ for every $e\in E$ and $\nrm{\pi(b,e)}\ge
\nrm{e}-\eta\nrm{b}$ for every $(e,b)\in E\times B$. \qed
\end{enumerate}

\fprop
It follows from (b) above that the composition of two $\eta$-admissible embeddings is also $\eta$-admissible.
Although we are not going to used them directly, two metric properties of $\eta$-admissible embeddings
crucial for the Radon-Nikodym and Schur properties of the Bourgain-Pisier quotient $Y_X/X$.
\begin{lemin}\label{kjdfklsdljfdw}
Suppose that the diagram $(\De)$  above is $\eta$-admissible. Then,
\begin{enumerate}
\item[(a)]  $\nrm{\overline{u}(b)}=\nrm{(b,0)+\ker(\pi)}=\inf_{s\in
S}\nrm{b+s}+\nrm{u(s)}$ for every $b\in B$. Consequently, $\nrm{\overline{u}}\le 1$, and if there is $\de\le 1$ such that $\nrm{u(s)}\ge \de\nrm{s}$  for every $s\in S$, then
$\nrm{\overline{u}(b)}\ge \de\nrm{b}$ for every $b\in B$. In other words, if $u$ is an isomorphic embedding then so is $\bar{u}$  with better isomorphic constant.
\item[(b)] Let $q:E_1\to E_1/j(E)$ be the natural quotient map. Suppose that $x_0,\dots,x_n\in E_1$ are
such that $x_0+\dots+x_n\in j(E)$. Then
$$\sum_{i=0}^n \nrm{x_i}\ge \nrm{\sum_{i=0}^n x_i}+(1-\eta)\sum_{i=0}^n \nrm{q(x_i)}.$$
\end{enumerate}
\end{lemin}
The fact in (b) is taken from \cite{Do} and it has an equivalent probabilistic reformulation in \cite{BP}. It
is the key to prove the following.
\teor\cite[Theorem 1.6.]{BP}\label{ijisjfjdsss}
Suppose that $(E_n)_n$ is a sequence of finite dimensional spaces, and suppose that $j_n:E_n\to E_{n+1}$ is
an $\eta$-admissible embedding for each $n$. Then the inductive limit of $(E_n,j_n)_n$ has the Schur and the
Radon-Nikodym properties.
\fteor

\subsection{One step extension} We finish this section with the following result, somehow  stating  that an  appropriate composition of
$\eta$-admissible diagrams is again $\eta$-admissible.
\lema \label{i4jrijeerdthtyolo}
\label{khwe4iothjiogff} Suppose that
\begin{equation}
\begin{tikzpicture}[descr/.style={fill=white,inner sep=2pt}]
\matrix (m) [matrix of math nodes, row sep=2em, column sep=2em, text height=1.5ex, text depth=0.25ex]
{ B_0 & X_0 & &  & X_2 \\
& & S_1 & B_1 &  \\
S_0 & E & & & X_1
\\
};

\path[->,font=\normalsize]

(m-1-1) edge node[above] {$\bar{u_0} $} (m-1-2)

(m-3-1) edge node[below]  {$u_0$}  (m-3-2)

(m-3-2) edge node[right] {$j_0$} (m-1-2)

(m-1-2) edge node[above] {$j_2$} (m-1-5)

(m-3-2) edge node[below] {$j_1$} (m-3-5)

(m-3-5) edge node[right] {$j$} (m-1-5)

(m-2-3) edge node[right] {$u_2$} (m-1-2)

(m-2-4) edge node[left] {$\overline{u_2}$} (m-1-5)

(m-2-3) edge node[right] {$u_1$} (m-3-2)

(m-2-4) edge node[left] {$\overline{u_1}$} (m-3-5)

 ;

\path[right hook->]

(m-3-1) edge node[right=8pt] {$(\De_0)$}     (m-1-1)

(m-2-3) edge  node[above=8pt] {$(\De.2)$} node[below=8pt] {$(\De.1)$}  (m-2-4)

;
\end{tikzpicture}
\end{equation}
is a commutative diagram such that:
\begin{enumerate}
\item[(1)] $(\De.0)$, $(\De.1)$ and $(\De.2)$ are $\eta$-admissible diagrams.
\item[(2)] $j:X_1\to X_2$ is an isometry.
\end{enumerate}
Then the diagram

\begin{equation}
\begin{tikzpicture}[descr/.style={fill=white,inner sep=2pt}]
\matrix (m) [matrix of math nodes, row sep=3em, column sep=3em, text height=1.5ex, text depth=0.25ex]
{ B_0 & X_2 \\
S_0 & X_1
\\
};

\path[->,font=\normalsize]

(m-1-1) edge node[above] {$j_2\circ\bar{u_0} $} (m-1-2)

(m-2-1) edge node[below]  {$j_1\circ u_0$}  (m-2-2)

(m-2-2) edge node[right] {$j$} (m-1-2)

 ;

\path[right hook->]

(m-2-1) edge
   (m-1-1)

;

\end{tikzpicture}
\end{equation}
is $\eta$-admissible.
\flema

\prue
Let $\pi:B_0\oplus_1 X_1\to X_2$, $\pi(b_0,x):=j_2(\bar{u_0}(b_0))+j(x)$, and for $i=0,1,2$, let  $\pi_i:
B_i\otimes E_i\to X_i$ be defined by $\pi_i(b,e):= \bar{u_i}(b)+j_i(e)$, where $E_0=E_1=E$, $E_2=X_0$ and
$B_2=B_1$. We have to check that $(\al.1)$, $(\al.2)$ and $(\al.3)$ in Proposition \ref{jiejrijeijifjdd} (a)
hold. By hypothesis $j$ is isometry and clearly $\nrm{j_1\circ u_0}=\nrm{u_0}\le \eta$, so we get $(\al.1)$.
\clam\label{nerhoifdohidgf}
$\pi(b_0,\pi_1(b_1,e))=\pi_2(b_1,\pi_0(b_0,e))$ for every $b_0\in B_0$, $b_1\in B_1$ and $e\in E$.
\fclam
\prucl
\begin{align*}
 \pi(b_0,\pi_1(b_1,e))=& \pi_2(b_1,\pi_0(b_0,e))=j_2(\bar{u_0}(b_0))+j(\pi(b_1,e))=j_2(\bar{u_0}(b_0))+j(\bar{u_1}(b_1)+ j_1(e))=\\
 =& j_2(\bar{u_0}(b_0)+j_0(e))+j(\bar{u_1}(b_1))= j_2(\bar{u_0}(b_0)+j_0(e))+ \bar{u_2}(b_1)=\\
 =& \pi_2(b_1,\bar{u_0}(b_0)+j_0(e))=\pi_2(b_1,\pi_0(b_0,e)).
 \end{align*}
\fprucl
It follows from this  that $\pi$ is onto.
\clam\label{hieroioijdfgdf}
$\ker (\pi)= \conj{(s_0,-j_1(u_0(s_0)))}{s_0\in S_0}=\conj{(b_0,\pi_1(0,-u_0(b_0)))}{b_0\in S_0}$.
\fclam
\prucl
The last equality follows from the fact that by definition, $\pi_1(0,-u_0(b_0))=-j_1(u_0(b_0))$. We prove now
the first equality. Fix  $s_0\in S_0$, and we work to prove that $\pi(s_0,-j_1(u_0(s_0)))=0$. Using the
commutativity of the diagram we obtain
\begin{align*}
\pi(s_0,-j_1(u_0(s_0)))= & j_2(\overline{u_0}(s_0))-j(j_1(u_0(s_0)))=j_2(\overline{u_0}(s_0))-j_2(j_0(u_0(s_0)))= \\
=& j_2(\overline{u_0}(s_0)-j_0(u_0(s_0)) )= j_2(0)=0.
\end{align*}
Now suppose that $\pi(b_0,g)=0$. Let $(b_1,e)\in B_1\times E$ be such that $\pi_1(b_1,e)=g$. Then, by Claim
\ref{nerhoifdohidgf}, it follows that
\begin{equation}
(b_1,\pi_0(b_0,e))\in \ker (\pi_2).
\end{equation}
And hence, $b_1\in S_1$ and $\pi_0(b_0,e)=-u_2(b_1)$.  It follows that
\begin{align*}
0= & \bar{u_0}(b_0)+j_0(e)+u_2(b_1)=\bar{u_0}(b_0)+j_0(e)+j_0(u_1(b_1))= \pi_0(b_0,e+u_1(b_1))
\end{align*}
So, $b_0\in S_0$ and $e+u_1(b_1)=-u_0(b_0)$. By applying $j_1$ to the last equality, we obtain that
\begin{align*}
g:=j_1(e)+\bar{u_1}(b_1)= -j_1(u_0(b_0)),
\end{align*}
as desired.
\fprucl
It follows readily that $(\al.3)$ holds. It rests to prove the property $(\al.2)$.

\clam\label{jioiojoijw34ef}
$\nrm{\pi(b_0,g)}=\inf_{s_0\in S_0}\nrm{b_0+s_0}+\nrm{g-j_1(u_0(s_0) )}$.
\fclam
\prucl
Fix $(b_1,e)\in B_1\times E$ such that $\pi_1(b_1,e)=g$. Then, by Claim \ref{nerhoifdohidgf} it follows that
$\pi(b_0,g)=\pi_2(b_1,\pi_0(b_0,e))$. Hence,
\begin{align*}
\nrm{\pi(b_0,g)}=& \nrm{\pi_2(b_1,\pi_0(b_0,e))}=\nrm{(b_1,\pi_0(b_0,e))+\ker(\pi_2)}=\\
=& \inf_{s_1\in S_1}\left( \nrm{b_1+s_1}+\nrm{\pi_0(b_0,e)-u_2(s_1)}\right)=\\
= & \inf_{s_1\in S_1} \left(\nrm{b_1+s_1}+\nrm{\pi_0(b_0,e)-j_0(u_1(s_1))}\right)=  \\
= & \inf_{s_1\in S_1} \left(\nrm{b_1+s_1}+\nrm{\pi_0(b_0,e-u_1(s_1))}\right)=  \\
=&  \inf_{s_1\in S_1}\left( \nrm{b_1+s_1}+\inf_{s_0\in S_0}\left( \nrm{b_0+s_0}+\nrm{e-u_1(s_1)-u_0(s_0)}\right)\right)= \\
=& \inf_{s_0\in S_0} \left( \nrm{b_0+s_0}+\inf_{s_1\in S_1}\left(\nrm{s_1+b_1}+\nrm{e-u_1(s_1)-u_0(s_0)} \right)\right)=\\
=& \inf_{s_0\in S_0} \left( \nrm{b_0+s_0}+\nrm{ (b_1,e-u_0(s_0))+\ker(\pi_1)} \right)=\\
=& \inf_{s_0\in S_0} \left( \nrm{b_0+s_0}+\nrm{\pi_1(b_1,e)+\pi_1(0,-u_0(s_0))} \right)=\\
=& \inf_{s_0\in S_0} \left( \nrm{b_0+s_0}+\nrm{g-j_1(u_0(s_0))} \right).
\end{align*}
\fprucl
From this we prove that $\pi$ is a metric surjection:  Fix $(b_0,g)\in B_0\times G$, and let $(b_1,e)\in
B_1\times E$ be such that $g=\pi_1(b_1,e)$. Then by the Claim \ref{hieroioijdfgdf}, it follows that
\begin{align*}
\nrm{(b_0,g)+\ker(\pi)}=& \nrm{(b_0,\pi_1(b_1,e))+\ker(\pi)}=\\
=& \inf_{s_0\in S_0}\left(\nrm{b_0+s_0}+\nrm{\pi_1(b_1,e)+\pi_1(0,-u_0(s_0))}\right)=\\
=& \inf_{s_0\in S_0}\left(\nrm{b_0+s_0}+\nrm{\pi_1(b_1,e-u_0(s_0))}\right)=\\
=& \inf_{s_0\in S_0}\left(\nrm{b_0+s_0}+\inf_{s_1\in S_1}\left(\nrm{b_1+s_1}+\nrm{e-u_0(s_0)-u_1(s_1)}\right)\right)=\\
=& \inf_{s_0\in S_0} \left( \nrm{b_0+s_0}+\nrm{g-j_1(u_0(s_0))} \right)=\nrm{\pi(b_0,g)},
\end{align*}
the last equality by Claim \ref{jioiojoijw34ef}.
\fprue

\section{The main result}\label{dkjferjeijr}
Our goal is to isometrically embed a given Banach space, not necessarily separable, into a $\mc
L_\infty$-space in such a way that the corresponding quotient has the Schur and the Radon-Nikodym properties.
Extending the approach of Bourgain and Pisier, we will find the $\mc L_\infty$-space as a direct, not
necessarily linear, limit of $\eta$-admissible embeddigs. The following is our main result.

\teor
\label{maintheorem} Every infinite dimensional Banach space $X$  can be isometrically embedded into a $\mc
L_\infty$-space $Y$ of the same density that $X$   such that the quotient $Y/X$ has the Radon-Nikodym and the
Schur properties.
\fteor

\cor
For every infinite cardinal number $\ka$ there is a $\mc L_\infty$-space of density $\ka$ with the
Radon-Nikodym and the Schur properties.
\fcor
\prue
For a fixed infinite cardinal number $\ka$, apply the Theorem \ref{maintheorem} to $X=\ell_1(\ka)$. Then the
corresponding superspace $Y$ is the desired space, since the required properties are three space properties.
\fprue

For the proof of Theorem \ref{maintheorem} we need the following two concepts.
\defi
Recall that the \emph{anti-lexicographical} ordering  $\prec$  on the family  $[\ka]^{<\om}$ of finite
subsets of $\ka$ is defined recursively as follows: $\buit \prec s$ for every non empty $s$, and
\begin{center}
$t\prec s $  if and only if   $\left\{ \begin{array}{ll}
\max t<\max s &  \text{ or}\\
\max t=\max s &\text{and }  t\setminus \{\max t\}\prec s \setminus \{\max s\}
\end{array}\right.$
\end{center}
\fdefi
  This is a well-ordering on $[\ka]^{<\om}$ that extends the inclusion relation $\cones$. We introduce some notation: For each
$\buit\cones t\con s$, we denote by $\bar{t}^{(s)}$ the immediate $\prec$-predecessor of $t$ in the family
$\mc P(s)$ of subsets of $s$, i.e.
$$\bar{t}^{(s)}:=\max_\prec\conj{u\con s}{u\prec t}.$$
Obviously this is well defined since $\mc P(s)$ is finite.  We write $\bar{t}$ to denote $\bar t^{(t)}$.

\defi
Recall that  a  \emph{directed system}   is   $((X_i)_{i\in I}, (j_{i_0,i_1})_{i_0\le _{I}i_1})$, where $X_i$
are Banach spaces, $<_I$ is a directed partial ordering, $j_{i_0,i_1}:X_{i_0}\to X_{i_1}$ are isometrical
embeddings, such that if $i_0\le_I i_i\le_I i_2$, then $j_{i_0,i_2}=j_{i_1,i_2}\circ j_{i_0,i_1}$, and such
that $j_{i,i}=\id_{X_i}$.

\fdefi

From now on we fix an infinite dimensional Banach space  $X$  of density $\ka$, and a  dense subset
$D=\conj{d_\al}{\al<\ka}$ of it. For each $s\in [\ka]^{<\om}$, let $X_s$ be the linear span of
$\{d_\al\}_{\al\in s}$. Fix also $\la>1$ and  $\eta<1$ such that $\la \cdot\eta < 1$.
\lema
\label{ni43hjoit4hjt} There is a direct systems $((E_s)_{s\in [\ka]^{<\om}},(j_{s,t})_{s\con t,\, s,t\in
[\ka]^{<\om}})$  and $(G_s)_{s\in [\ka]^{<\om}}$   such that:
\begin{enumerate}
\item[(1)] $G_s\con E_s$ are Banach spaces, $E_\buit=X$.
\item[(2)] Each $j_{s,t}:E_s\to E_t$ is an   $\eta$-admissible isometrical embedding such that $j_{s,t}E_s$
has finite codimension in $E_t$.
\item[(3)] $G_s$ is $\la$-isomorphic to $\ell_\infty^{\dim G_s}$.
\item[(4)] $\bigcup_{t\cones s}j_{t,s}(G_t)\cup j_{\buit,s}(X_s)\con G_s$.
\end{enumerate}
\flema
We are ready now to give a proof   of Theorem \ref{maintheorem} from this lemma.

\prue[Proof of Theorem \ref{maintheorem}]
Fix  $((E_s)_{s\in [\ka]^{<\om}},(j_{s,t})_{s\con t,\, s,t\in [\ka]^{<\om}})$  and $(G_s)_{s\in
[\ka]^{<\om}}$  as in Lemma \ref{ni43hjoit4hjt}. Let  $E$ be the completion of the inductive limit of
$((E_s)_{s\in [\ka]^{<\om}}, (j_{t,s})_{t\con s,\, t,s\in [\ka]^{<\om}})$. Because of property (4) in Lemma
\ref{ni43hjoit4hjt}, it follows that $((G_s)_{s\in \mathcal F},(j_{t,s}\rest G_t))_{t\con s\in [\ka]^{<\om}}$
is also a directed system of finite dimensional normed spaces $G_s$ which are $\la$-isomorphic to
$\ell_\infty^{\dim G_s}$. Let $Y$ be the completion of the corresponding direct limit $\lim_{s\in
[\ka]^{<\om}}G_s$. It is clear that $Y$ can be isometrically imbedded into $E$, while there is a natural
isometric embedding of $X$ into $Y$:   $X$ is the completion of the direct limit $((X_s)_{s\in
[\ka]^{<\om}},(i_{t,s})_{t\con s\in [\ka]^{<\om}})$, where $i_{t,s}:X_t\to X_s$ is the inclusion map. For
each finite subset $s$ of $\ka$, let $g_{s}:X_s\to G_s$ be $g_s:=j_{\buit,s}\rest X_s$, which is well defined
by (4). This is obviously an isometric embedding such that $j_{t,s}\circ g_{s}=g_s\circ i_{t,s}$ for every
$t\con s$, and hence $X$ isometrically embeds into $Y$.

If we denote by $j_{s,\infty}: G_s \to Y$ the corresponding limit of $(j_{s,t})_{s\con t}$, then
$\bigcup_{s\in [\ka]^{<\om}}j_{s,\infty} (G_s)$ is dense in $Y$. It follows that $Y$ is a $\mc
L_{\infty,\la}$-space. Since each $G_s$ is finite dimensional, it follows that $Y$ has density at most
$|[\ka]^{<\om}|=\ka$. Since $X$ isometrically embeds into $Y$, the density of $Y$ has to be $\ka$.   Let us
see that $Y/X$ has the Radon-Nikodym and the Schur  properties: We use that $Y/X$ is naturally isometrically
embedded into $E/X$, and we prove that $E/X$ has these two properties. Observe that these two properties are
properties of separable subspaces of $E/X$. So let $Z\con E/X$ be a separable subspace of $E/X$. By
construction, we can find a sequence $(s_n)_{n\in \N}$ of elements of $\mathcal F$ such that $s_n\con
s_{n+1}$, and such that $Z$ is a subspace of the closure of the quotient
$$\left(\lim_{n\to \infty} ((E_{s_n})_{n\in
\N},(j_{s_n s_{n+1}})_{n\in \N})\right)/X.$$ This quotient can be  naturally isometrically identified with
the inductive limit of finite dimensional spaces $((E_{s_n}/j_{\buit, s_n}X )_{n\in
\N},(\overline{j_{s_n,s_{n+1}}})_{n\in \N})$, which, by Theorem \ref{ijisjfjdsss},  has the two required
properties.
\fprue
The existence of the direct system  in Lemma \ref{ni43hjoit4hjt} is based on the following local
construction.
\lema\label{j4irjtiojgghff}
For every finite subset  $s$ of $\ka$  there are
\begin{equation}
\label{kjhuhurt}\text{$(E_t^{(s)})_{t\con s}$, $(G_t^{(s)})_{t\con s}$,  $(j_{u,t}^{(s)})_{u\prec
t,\,u,t\con s}$ and $(k_{u}^{(t,s)})_{u\con t\con s}$}
\end{equation}
such that
\begin{enumerate}
\item[(A)] (Local directed system)  For every finite subset $s$ of $\ka$  one has that
$$((E_t^{(s)})_{t\con s}, (j_{u,t}^{(s)})_{u\prec t, \, u,t\con s})$$
  is a (finite)   system of $\eta$-admissible isometrical embeddings such that $j_{u,t}^{(s)}E_u^{(s)}$
  has finite codimension in $E_t^{(s)}$.

\item[(B)] (Transition directed system) For every  finite subset $u$ of $\ka$   one has that
 $$((E_u^{(s)})_{u\con s\in [\ka]^{<\om}}, (k_{u}^{(t,s)})_{u\con t\con s\in [\ka]^{<\om}})$$
 is a  system of $\eta$-admissible isometrical embeddings such that $k_u^{(t,s)}E_u^{(t)}$ has finite
 codimension in $E_u^{(s)}$.
\item[(C)] (Coherence property) For every $v\con u \con t\con s$ one has that
\begin{equation}
\label{erjijgijgfj}  k_{u}^{(t,s)}\circ j_{v,u}^{(t)}= j_{v,u}^{(s)}\circ k_{v}^{(t,s)}.
\end{equation}
\item[(D)] For every $t\con s\in [\ka]^{<\om}$ one has that $G_t^{(s)}\con E_t^{(s)}$ is finite dimensional and
\begin{equation}
\label{cassegffb} d(G_t^{(s)},\ell_\infty^{\dim G_t^{(s)}})\le \la.
\end{equation}
\item[(E)] For every $u\con t\con s$ one has that
\begin{equation}
\label{ljuhbbb} k_{u}^{(t,s)} (G_{u}^{(t)})=G_{u}^{(s)}.
\end{equation}
\item[(F)] For every $\buit \cones t\con s$  one has that
\begin{equation}
\label{rgjrjghrhhhhff} j_{\bar t,t}^{(s)} (G_{\bar t}^{(s)} )\cup j_{\buit,t}^{(s)} (X_{t})\con G_t^{(s)}.
\end{equation}
\end{enumerate}
\flema
We postpone its proof, and we pass to prove Lemma \ref{ni43hjoit4hjt}.
\prue[Proof of Lemma \ref{ni43hjoit4hjt}]
Fix $(E_t^{(s)})_{t\con s}$, $(G_t^{(s)})_{t\con s}$,  $(j_{u,t}^{(s)})_{u\prec t,\,u,t\con s}$ and
$(k_{u}^{(t,s)})_{u\con t\con s}$ as in   Lemma \ref{j4irjtiojgghff}.  For each finite subset $s$  of $\ka$
set  $E_s:=E_s^{(s)}$ and $G_s:=G_s^{(s)}$. Given $t\con s$  let $j_{t,s}: E_t\to E_s$   be
\begin{equation} \label{ueuuurjrihjh}
j_{t,s}:=  j_{t,s}^{(s)}\circ k_{t}^{(t,s)}.
\end{equation}
It follows from properties (A) and (B) in Lemma \ref{j4irjtiojgghff} that $j_{t,s}$ is an $\eta$-admissible
embedding such that $j_{t,s}E_t$ has finite codimension in $E_s$.
\clam
 $((E_{s})_{s\in [\ka]^{<\om}},(j_{t,s})_{t\con s\in [\ka]^{<\om}})$ is a directed system of $\eta$-admissible isometrical
 embeddings.
\fclam

\prucl
Suppose that $u\con t\con s$. Then
\begin{align*}
j_{t,s}\circ j_{u,t}= & j_{t,s}^{(s)}\circ( k_{t}^{(t,s)}\circ j_{u,t}^{(t)})\circ k_{u}^{(u,t)}=_{\eqref{erjijgijgfj}}
j_{t,s}^{(s)}\circ( j_{u,t}^{(s)}\circ  k_{u}^{(t,s)})\circ   k_{u}^{(u,t)}= \\
= & (j_{t,s}^{(s)}\circ j_{u,t}^{(s)})\circ  (k_{u}^{(t,s)}\circ   k_{u}^{(u,t)})=
j_{u,s}^{(s)}\circ k_{u}^{(u,s)} =j_{u,s}
\end{align*}
\fprucl
For each $s\in [\ka]^{<\om}$, let $G_s=G_s^{(s)}$.
\clam
\begin{enumerate}
\item[(a)] $d(G_s,\ell_\infty^{\dim G_s})\le \la$, i.e. $G_s$ is $\la$-isomorphic to $\ell_\infty^{\dim
G_s}$.
\item[(b)]  For every $t\con s$  one has that
$j_{t,s}G_t\con G_s$.
\end{enumerate}
\fclam
\prucl

(a) follows from \eqref{cassegffb} in (D). (b): By \eqref{ljuhbbb} and \eqref{rgjrjghrhhhhff}  one has that
$$j_{t,s} (G_t)=
j_{t,s}^{(s)}\circ k_{t}^{(t,s)} (G_t^{(t)})= j_{t,s}^{(s)} G_t^{(s)}\con G_s^{(s)}=G_s.$$
 \fprucl
\fprue
It only rests to give a proof of Lemma  \ref{j4irjtiojgghff}.
\prue[Proof of Lemma \ref{j4irjtiojgghff}]

Fix $s\in [\ka]^{<\om}$.  We define $\pe$-recursively on $s$ all the objects in \eqref{kjhuhurt} together
with  an integer $n_t\in \N$, $S_t\con \ell_\infty^{n_t}$ and $\upsilon_t:S_t\to E_{\bar t}^{(t)}$  for each
 $\buit\cones t\con s$  such that

\begin{enumerate}
\item[(a)] $E_\buit^{(s)}=X$, $k_{\buit}^{(t,s)}=\mr{Id}_{X}$.
\item[(b)]  $\upsilon_t:S_t\to E_{\bar t}^{(t)}$ is an isomorphism with $\nrm{\upsilon_t}\le \eta$, $\nrm{\upsilon_t^{-1}}\le \la$, and
\begin{equation}
\label{jutrihgf} \upsilon_t(S_t)=\langle j_{u,\bar t}^{(t)}(G_u^{(t)})\cup j_{\buit,\bar{t}}^{(t)}(X_t)\rangle\con E_{\bar t}^{(t)},
\end{equation}
where $u=\overline{\overline{t}}^{(t)}$ is the $\prec$-penultimate element in $\mc P(t)$, if $|t|>1$ and
$u=\buit$, if $|t|=1$ (and hence $\bar t=\buit$).
\item[(c)]  $ E_t^{(s)}=(\ell_\infty^{n_t}\oplus_1 E_{\overline{t}^{(s)}}^{(s)})/N_{v_t^{(s)}}$, the diagram
\begin{equation*}
\begin{tikzpicture}[descr/.style={fill=white,inner sep=2pt}]
\matrix (m) [matrix of math nodes, row sep=3em, column sep=2.5em, text height=1.5ex, text depth=0.25ex]
{\ell_\infty^{n_t} & E_t^{(s)}=(\ell_\infty^{n_t}\oplus_1 E_{\overline{t}^{(s)}}^{(s)})/N_{v_t^{(s)}} \\
S_t & E_{\overline{t}^{(s)}}^{(s)}
\\
E_{\overline{t}}^{(t)} & E_{\overline{t}}^{(s)} \\
};

\path[->,font=\normalsize]

(m-1-1) edge node[above] {$\overline{v_t^{(s)}} $} (m-1-2)

(m-2-1) edge node[below]  {$v_t^{(s)}$}  (m-2-2)

        edge node[left] {$v_t$}  (m-3-1)

(m-2-2) edge node[right] {$j_{\overline{t}^{(s)},t}^{(s)}$} (m-1-2)

(m-3-1) edge node[below] {$k_{\overline{t}}^{(t,s)}$} (m-3-2)

(m-3-2) edge node[right] {$j_{\overline{t},\overline{t}^{(s)}}^{(s)}$} (m-2-2)

 ;

\path[right hook->]

(m-2-1) edge node[right=45pt] {$(\De)$} (m-1-1)

;

\end{tikzpicture}
\end{equation*}
is commutative,$(\De)$ is a canonical $\eta$-admissible diagram, and
\begin{equation}
G_{t}^{(s)}= \overline{v_t^{(s)}}(\ell_\infty^{n_t})=\conj{(x,0)+N_{t}^{(s)}}{x\in \ell_\infty^{n_t}}.\label{rjtirjirjr3}
\end{equation}

\item[(d)]  For every $u\prec t$ subset of $s$, we have that
\begin{equation}
\label{rtkjritjjhjjgg} j_{u,t}^{(s)}=j_{\overline{t}^{(s)},t}^{(s)}\circ j_{\overline{t}^{(s)},t}^{(s)}.
\end{equation}
\item[(e)] For every $u\con t\con s$, $k_{u}^{(t,s)}:E_{u}^{(t)}\to E_{u}^{(s)}$ satisfies that for  every $(x,y)\in
\ell_\infty^{n_u}\times E_{\bar{u}^{(t)}}^{(t)}$ by
\begin{equation}
\label{lejritjgjg}
k_u^{(t,s)}((x,y)+N_u^{(t)})=(x,j_{\bar{u}^{(t)},\bar{u}^{(s)}}^{(s)}\circ k_{\bar u^{(t)}}^{(t,s)}(y))+N_u^{(s)}.
\end{equation}
The requirement in (e) can be fulfilled   because the commutativity of the following diagram:
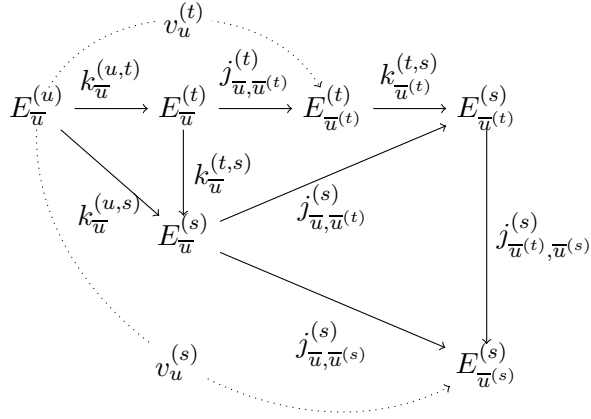
\begin{equation*}
\begin{tikzpicture}[descr/.style={fill=white,inner sep=2pt}]
\matrix (m) [matrix of math nodes, row sep=3em, column sep=2.5em, text height=1.5ex, text depth=0.25ex]
{E_{\overline{u}}^{(u)}& E_{\overline{u}}^{(t)} & E_{\overline{u}^{(t)}}^{(t)} &  E_{\overline{u}^{(t)}}^{(s)} \\
&  E_{\overline{u}}^{(s)} & & \\
& & &  E_{\overline{u}^{(s)}}^{(s)}
\\ };
\path[->,font=\normalsize]

(m-1-1) edge node[auto] {$ k_{\overline{u}}^{(u,t)}$} (m-1-2)

(m-1-2) edge node[auto] {$j_{\overline{u},\overline{u}^{(t)}}^{(t)}$} (m-1-3)

 (m-1-3) edge node[auto] {$ k_{\overline{u}^{(t)}}^{(t,s)} $} (m-1-4)

  (m-1-4) edge node[auto] {$j_{\overline{u}^{(t)},\overline{u}^{(s)}}^{(s)}$} (m-3-4)

 (m-1-1) edge node[ below=4pt ] {$k_{\overline{u}}^{(u,s)}$} (m-2-2)

 (m-2-2) edge node[below=4pt] {$j_{\overline{u},\overline{u}^{(s)}}^{(s)}$} (m-3-4)

  (m-2-2) edge node[below] {$j_{\overline{u},\overline{u}^{(t)}}^{(s)}$} (m-1-4)

  (m-1-2) edge node[auto] {$k_{\overline{u}}^{(t,s)}$} (m-2-2) ;

\path[->,dotted]

(m-1-1) edge[bend right=60] node[descr]{$v_{u}^{(s)}$} (m-3-4)

(m-1-1) edge [bend left=60] node[descr] {$v_{u}^{(t)}$} (m-1-3) ;
\end{tikzpicture}
\end{equation*}
\end{enumerate}
It rests to check that the conditions (A)--(F) hold:

\noindent (A): It is clear from the definition of $j_{\bar{t}^{(s)},t}^{(s)}$ is an $\eta$-admissible
isometrical embedding, and it follows from the equality in  \eqref{rtkjritjjhjjgg} that
$j_{u,t}^{(s)}=j_{v,t}^{(s)}\circ j_{u,v}^{(s)}$ for every $u,v\con s$ with $u\prec v \prec t$.

\noindent (C): Let $v\cones u \con t\con s$. We want to prove that $ k_{u}^{(t,s)}\circ j_{v,u}^{(t)}=
j_{v,u}^{(s)}\circ k_{v}^{(t,s)}$. Using that, by inductive hypothesis, the left side of the following
diagram is commutative,
$$
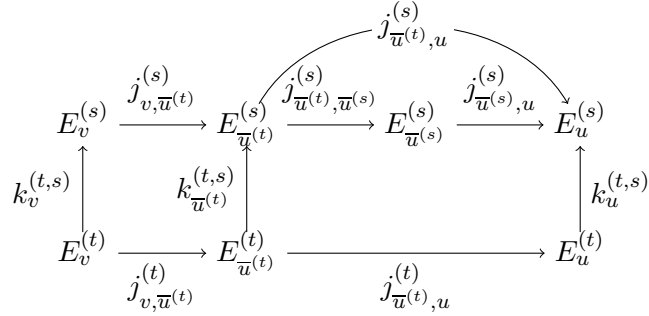
\begin{tikzpicture}[descr/.style={fill=white,inner sep=2pt}]
\matrix (m) [matrix of math nodes, row sep=3em, column sep=3em, text height=1.5ex, text depth=0.25ex]
{E_v^{(s)} & E_{\overline{u}^{(t)}}^{(s)} & E_{\overline{u}^{(s)}}^{(s)} & E_u^{(s)}\\
E_v^{(t)} & E_{\overline{u}^{(t)}}^{(t)} & & E_{u}^{(t)}
\\
};

\path[->,font=\normalsize]

(m-1-1) edge node[above] {$j_{v,\overline{u}^{(t)}}^{(s)}$} (m-1-2)

(m-1-2) edge node[above] {$j_{\overline{u}^{(t)},\overline{u}^{(s)}}^{(s)}$} (m-1-3)

(m-1-3) edge node[above] {$j_{\overline{u}^{(s)},u}^{(s)}$} (m-1-4)

(m-1-2) edge[bend left=60] node[descr] {$j_{\overline{u}^{(t)},u}^{(s)}$}   (m-1-4)

(m-2-1) edge node[left] {$k_v^{(t,s)}$}  (m-1-1)

(m-2-2) edge node[left] {$k_{\overline{u}^{(t)}}^{(t,s)}$} (m-1-2)

(m-2-1) edge node[below] {$j_{v,\overline{u}^{(t)}}^{(t)}$}  (m-2-2)

(m-2-2) edge node[below] {$j_{\overline{u}^{(t)},u}^{(t)}$}   (m-2-4)

(m-2-4) edge node[right] {$k_u^{(t,s)}$} (m-1-4)

 ;
\end{tikzpicture}
$$
it suffices to prove that $k_{u}^{(t,s)}\circ j_{\bar{u}^{(t)},u}^{(t)}=j_{\bar{u}^{(t)},u}^{(s)}\circ
k_{\bar{u}^{(t)}}^{(t,s)}$. Let $x\in E_{\bar u^{(t)}}^{(t)}$. Then by (e),
\begin{align*}
k_{u}^{(t,s)}\circ j_{\bar{u}^{(t)},u}^{(t)}(x)= & k_u^{(t,s)}((0,x)+N_u^{(t)})=(0,j_{\bar u^{(t)},\bar{u}^{(s)}}^{(s)}\circ
k_{\bar u^{(t)}}^{(t,s)}(x))+N_u^{(s)}\\
j_{\bar{u}^{(t)},u}^{(s)}\circ
k_{\bar{u}^{(t)}}^{(t,s)}(x)=&  j_{\bar{u}^{(s)},u}^{(s)} \circ j_{\bar{u}^{(t)},\bar{u}^{(s)}}^{(s)}\circ
k_{\bar{u}^{(t)}}^{(t,s)}(x)=  (0,j_{\bar{u}^{(t)},\bar{u}^{(s)}}^{(s)}\circ
k_{\bar{u}^{(t)}}^{(t,s)}(x) )+ N_{u}^{(s)}.
\end{align*}

\noindent (B): Suppose that $v\con u\con t\con s$. We have to see that $k_{v}^{(u,s)}= k_{v}^{(t,s)}\circ
k_{v}^{(u,t)}$.    Recall that from (e) it follows that  for every $(x,y)\in \ell_\infty^{n_v}\times E_{\bar
v^{(u)}}$ and for every $(x,z)\in \ell_\infty^{n_v}\times E_{\bar v^{(t)}}$ one has that
\begin{align}
\label{rtjrtjgijgjgg} k_v^{(u,t)}((x,y)+N_v^{(u)})=&(x, j_{\bar v^{(u)},\bar v^{(t)}}^{(t)}\circ k_{\bar v^{(u)}}^{(u,t)}(y))+N_v^{(t)}
\\
k_v^{(u,s)}((x,y)+N_v^{(u)})=&(x, j_{\bar v^{(u)},\bar v^{(s)}}^{(s)}\circ k_{\bar v^{(u)}}^{(u,s)}(y))+N_v^{(s)} \\
k_v^{(t,s)}((x,z)+N_v^{(t)})=&(x, j_{\bar v^{(t)},\bar v^{(s)}}^{(s)}\circ k_{\bar v^{(t)}}^{(t,s)}(z))+N_v^{(s)}.
\end{align}
Hence, using inductively  (C),
\begin{align*}
\label{rtjrtjgijgjgg} k_v^{(t,s)}\circ k_v^{(u,t)}((x,y)+N_v^{(u)})=&(x, j_{\bar v^{(t)},\bar v^{(s)}}^{(s)}
\circ k_{\bar v^{(t)}}^{(t,s)}\circ j_{\bar v^{(u)},\bar v^{(t)}}^{(t)}\circ k_{\bar v^{(u)}}^{(u,t)}(y))+N_v^{(t)}=
\\
= & (x, j_{\bar v^{(t)},\bar v^{(s)}}^{(s)}
\circ j_{\bar v^{(u)},\bar v^{(t)}}^{(s)}\circ k_{\bar v^{(u)}}^{(t,s)}\circ  k_{\bar v^{(u)}}^{(u,t)}(y))+N_v^{(t)}   \\
= & (x, j_{\bar v^{(u)},\bar v^{(s)}}^{(s)}
 \circ k_{\bar v^{(u)}}^{(u,s)}(y))+N_v^{(t)}   \\
= & k_v^{(u,s)}((x,y)+N_v^{(u)}).\\
\end{align*}
We now prove that $k_u^{(t,s)}$ is an $\eta$-admissible isometrical embedding: By inductive hypothesis  the
composition $j:E_{\bar u^{(t)}}^{(t)}\to E_{\bar u^{(s)}}^{(s)}$,  $j:=j_{\bar u^{(t)},\bar{ u}^{(s)}}
^{(s)}\circ k_{\bar u^{(t)}}^{(t,s)}$, is an $\eta$-admissible isometrical embedding. We then fix $S\con B$,
 $\nu:S\to E_{\bar u^{(t)}}^{(t)}$  and $\overline{\nu}:B\to
E_{\overline{u}^{(s)}}^{(s)}$  such that

\begin{equation*}
\begin{tikzpicture}[descr/.style={fill=white,inner sep=2pt}]
\matrix (m) [matrix of math nodes, row sep=3em, column sep=3em, text height=1.5ex, text depth=0.25ex]
{ B &  E_{\bar u^{(s)}}^{(s)} \\
S & E_{\bar u^{(t)}}^{(t)}
\\
};

\path[->,font=\normalsize]

(m-1-1) edge node[above] {$\bar{\nu} $} (m-1-2)

(m-2-1) edge node[below]  {$\nu$}  (m-2-2)

(m-2-2) edge node[right] {$j=j_{\bar u^{(t)},\bar{ u}^{(s)}} ^{(s)}\circ k_{\bar u^{(t)}}^{(t,s)}$} (m-1-2)

 ;

\path[right hook->]

(m-2-1) edge node[right=13pt] {$(\De_0)$}   (m-1-1)

;

\end{tikzpicture}
\end{equation*}
is an  $\eta$-admissible diagram.  It follows  by (C), that the following diagram is commutative:

\begin{equation*}
\begin{tikzpicture}[descr/.style={fill=white,inner sep=2pt}]
\matrix (m) [matrix of math nodes, row sep=2em, column sep=2em, text height=1.5ex, text depth=0.25ex]
{ B & E_{\bar{u}^{(s)}}^{(s)} & &  & E_u^{(s)}=(\ell_\infty^{n_u}\oplus_1 E_{\bar u^{(s)}}^{(s)})/N_u^{(s)} \\
& & S_u & \ell_\infty^{n_u} &  \\
S & E_{\bar u^{(t)}}^{(t)} & & & E_u^{(t)}=(\ell_\infty^{n_u}\oplus_1 E_{\bar u^{(t)}}^{(t)})/N_u^{(t)}
\\
};
\path[->,font=\normalsize]

(m-1-1) edge node[above] {$\bar{u_0} $} (m-1-2)

(m-3-1) edge node[below]  {$u_0$}  (m-3-2)

(m-3-2) edge node[right] {$j$} (m-1-2)

(m-1-2) edge node[above] {$j_{\overline{u}^{(s)},u}^{(s)}$} (m-1-5)

(m-3-2) edge node[below] {$j_{\overline{u}^{(t)},u}^{(t)}$} (m-3-5)

(m-3-5) edge node[right] {$k_{u}^{(t,s)}$} (m-1-5)

(m-2-3) edge node[right=3pt] {$v_u^{(s)}$} (m-1-2)

(m-2-4) edge node[left=18pt] {$\overline{v_u^{(s)}}$} (m-1-5)

(m-2-3) edge node[right] {$v_u^{(t)}$} (m-3-2)

(m-2-4) edge node[right=2pt, above=1pt] {$\overline{v_u^{(t)}}$} (m-3-5)

 ;
\path[right hook->]

(m-3-1) edge node[right=8pt] {$(\De_0)$}     (m-1-1)

(m-2-3) edge  node[above=8pt] {$(\De.2)$} node[below=8pt] {$(\De.1)$}  (m-2-4) ;
\end{tikzpicture}
\end{equation*}
Since $(\De_0)$, $(\De_1)$ and $(\De_2)$ are $\eta$-admissible diagram, we conclude from   Lemma
\ref{khwe4iothjiogff} that $k_{u}^{(t,s)}$ is an $\eta$-admissible embedding.

\noindent (D) is clear by definition of $G_t^{(s)}$.

\noindent (E): Fix $u\con t\con s$.  Then
\begin{align*}
k_{u}^{(t,s)}(G_u^{(t)})=& \conj{ k_{u}^{(t,s)}((x,0)+N_{u}^{(t)})}{x\in \ell_\infty^{n_u}}=
\conj{ (x,0)+N_{u}^{(s)}}{x\in \ell_\infty^{n_u}}=G_{u}^{(s)}.
\end{align*}
\noindent (F): Let $t\con s$. We have to prove the inclusion in \eqref{rgjrjghrhhhhff}. Notice that the
diagram
\begin{equation*}
\begin{tikzpicture}[descr/.style={fill=white,inner sep=2pt}]
\matrix (m) [matrix of math nodes, row sep=3em, column sep=3em, text height=1.5ex, text depth=0.25ex]
{\ell_{\infty}^{n_t} & E_t^{(s)}=(\ell_\infty^{n_t}\oplus_1 E_{\overline{u}^{(s)}}^{(s)})/N_t^{(s)}\\
S_t & E_{\overline{u}^{(t)}}^{(s)}
\\
};

\path[->,font=\normalsize]

(m-1-1) edge node[above] {$\overline{v_t^{(s)}}$} (m-1-2)

(m-2-1) edge node[below]  {$v_t^{(s)}$}  (m-2-2)

(m-2-2) edge node[right] {$j_{\overline{t}^{(s)},t}^{(s)}$} (m-1-2)

 ;

\path[right hook->]

(m-2-1) edge  (m-1-1)

;

\end{tikzpicture}
\end{equation*}
is commutative. Let $u\con t$ be the immediate $\prec$-predecessor of $\overline{t}$ in $t$, if $|t|>1$, and
let $u=\buit$ otherwise. Then by (b),
\begin{align*}
G_t^{(s)}=& \conj{(x,0)+N_t^{(s)}}{s\in S_t}= \bar\upsilon_t^{(s)}(\ell_\infty^{n_t}) \supseteq \\
\supseteq & \bar\upsilon_t^{(s)}(S_t)=
j_{\overline{t}^{(s)},t}^{(s)}\circ \upsilon_t^{(s)}(S_t)= \\
= & j_{\overline{t}^{(s)},t}^{(s)}\circ j_{\bar{t},\bar{t}^{(s)}}^{(s)}\circ k_{\bar{t}}^{(t,s)}\circ \upsilon_t  (S_t)=
j_{\bar{t},t}^{(s)}\circ k_{\bar{t}}^{(t,s)}\circ \upsilon_t(S_t)= k_{t}^{(t,s)}\circ  j_{\bar{t},t}^{(t)}\circ \upsilon_t(S_t)=\\
= & k_{t}^{(t,s)}\circ  j_{\bar{t},t}^{(t)}\left\langle
j_{u,\bar{t}}^{(t)}(G_u^{(t)})\cup j_{\buit,\bar{t}}^{(t)}(X_t)\right\rangle =\\
=&  \left\langle k_{t}^{(t,s)}\circ  j_{\bar{t},t}^{(t)}\circ
j_{u,\bar{t}}^{(t)}(G_u^{(t)})\cup k_{t}^{(t,s)}\circ
  j_{\bar{t},t}^{(t)}\circ j_{\buit,\bar{t}}^{(t)}(X_t)\right\rangle= \\
=& \left\langle  k_{t}^{(t,s)}\circ
j_{u,t}^{(t)}(G_u^{(t)})\cup k_{t}^{(t,s)}\circ j_{\buit,t}^{(t)}(X_t)\right\rangle=\\
= &  \left\langle
j_{u,t}^{(s)}\circ k_{u}^{(t,s)}(G_u^{(t)})\cup  j_{\buit,t}^{(s)}\circ k_{\buit}^{(t,s)}(X_t)\right\rangle= \\
= &  \left\langle
j_{u,t}^{(s)}(G_u^{(s)})\cup  j_{\buit,t}^{(s)}(X_t)\right\rangle.
\end{align*}
\fprue

\end{document}